\documentclass[12pt]{amsart}
\usepackage{geometry} 
\usepackage[parfill]{parskip}    
\usepackage{graphicx}
\usepackage{amssymb}
\newcommand{\R}{\mathbb{R}}

\newtheorem{lma}{Lemma}

\newtheorem{thm}{Theorem}

\usepackage{amsfonts}

\begin{document}

\title[Saddle Point Method for Helmholtz Equations]{A Saddle Point Numerical Method for Helmholtz Equations}

\author{Russell B. Richins}
\address{Department of Mathematics and Computer Science, Thiel College, Greenville, Pennsylvania, 16125}
\email{rrichins@thiel.edu}


\maketitle

\begin{abstract} In a previous work, the author and D.C.\@ Dobson proposed a numerical method for solving the complex Helmholtz equation based on the minimization variational principles developed by Milton, Seppecher, and Bouchitt\'{e}.  This method results in a system of equations with a symmetric positive definite coefficient matrix, but at the same time requires solving simultaneously for the solution and its gradient.  Herein is presented a method based on the saddle point variational principles of Milton, Seppecher, and Bouchitt\'{e}, which produces symmetric positive definite systems of equations, but eliminates the necessity of solving for the gradient of the solution.  The result is a method for a wide class of Helmholtz problems based completely on the Conjugate Gradient algorithm.  
\end{abstract}

\section{Introduction} The Helmholtz equation
\[\nabla\cdot L\nabla u=Mu,\]
is useful in modeling wave propagation in problems arising from many different physical situations.  We will focus only on the homogeneous equation for simplicity and brevity, but the methods presented here can easily be extended to the non-homogeneous case.  Suppose we wish to solve the Helmholtz equation in a domain $\Omega\subset\R^d$, and assume that $L$ and $M$ are complex-valued functions.  A common source of numerical methods for solving this equation is the variational principle 
\begin{equation}\label{reg_var}\int_\Omega\left[-L\nabla u\cdot\nabla\bar v-Mu\bar v\right]dx=0\ \forall\ v\in H_0^1(\Omega).\end{equation}
Since this is a stationary principle, the resulting system of equations is often indefinite, and indefinite systems are generally more difficult to solve than a system of equations having a positive definite coefficient matrix.  

Because of the challenges in solving these indefinite systems, there has been much work devoted to solving the Helmholtz equation by replacing the indefinite systems with equivalent symmetric positive definite linear systems.  Classical examples of such approaches are the CGNR and CGNE methods \cite{Saad_2006}, based on solving normal equations associated with the original system.  While such approaches produce positive definite systems, the normal equations are often poorly conditioned and preconditioning can be difficult.  Another related approach is First Order System Least Squares (FOSLS) \cite{Cai_1994} \cite{Cai_1997} \cite{Lee_2000}, which converts the second order equation into an equivalent system of first order equations and then solves a least squares problem for this system.  The method presented here also produces positive definite systems of equations, but it does so without reformulation as a least squares problem.

When iterative methods are employed to solve a system of linear equations, it is usually necessary to precondition the original system in order to speed up convergence.  A great deal of  work has been dedicated to formulating effective preconditioning strategies for the linear systems resulting from discretizations  of the Helmholtz equation \cite{Gander_2012}.  One approach that has seen much success is the Shifted Laplacian preconditioner \cite{Erlangga_2004} \cite{Erlangga_2006} \cite{vanGizjen_2007}.  In this approach, the precoditioner for the system of equations corresponding to $\Delta u+k^2u=0$ is the matrix corresponding to the ``shifted" equation $\Delta u+(\alpha+i\beta)u=0$.  If the imaginary shift $\beta$ is large enough, multigrid methods are expected to be successful in solving the shifted problem, and if $\alpha\approx 1$, the shifted operator should be a good preconditioner for the original problem.  While this approach is often effective, in \cite{Vecharynski_2013} the authors point out the advantages in using a preconditioner that is symmetric positive-definite.  When the preconditioning matrix is not positive-defininte, the coefficient matrix of the preconditioned system is not symmetric with respect to any inner product, which limits the methods available for solving the resulting system.  The solution suggested in \cite{Vecharynski_2013} is to use an approximation of the absolute value of the original coefficient matrix as preconditioner.  In the method proposed here, both the matrices and the suggested preconditioners are symmetric positive-definte, and therefore a wide range of Krylov subspace methods is available.  In particular, we shall demonstrate the results obtained with Conjugate Gradient, which has a short recurrence and is very simple to implement and parallelize.

As a background to this approach, we start with \cite{Milton_2009}, where Milton, Seppecher, and Bouchitt\'{e} developed variational principles that apply to the Helmholtz equation above, as well as the time-harmonic Maxwell equations and the equations of linear elasticity in lossy materials.  To derive these variational principles, we first define the dual variable
\[v=iL\nabla u.\]
Then
\[\left(\begin{array}{cc}
L & 0 \\
0 & M \\
\end{array}\right)\left(\begin{array}{c}
\nabla u \\
u \\
\end{array}\right)=\left(\begin{array}{c}
L\nabla u \\
Mu \\
\end{array}\right)=\left(\begin{array}{c}
-iv \\
-i\nabla\cdot v \\
\end{array}\right),\]
or equivalently, 
  \[\mathcal{G}=Z\mathcal{F},\]
where 
\[\mathcal{F}=\left(\begin{array}{c}
\nabla u \\
u \\
\end{array}\right),\ \mathcal{G}=\left(\begin{array}{c}
-iv \\
-i\nabla\cdot v \\
\end{array}\right),\ Z=\left(\begin{array}{cc}
L & 0 \\
0 & M \\
\end{array}\right).\]
For a complex quantity $z$, we will write $z'=\mbox{Re}(z)$ and $z''=\mbox{Im}(z)$.  Taking real and imaginary parts, the constitutive relation becomes
\[\mathcal{G}'=Z'\mathcal{F}'-Z''\mathcal{F}''\mbox{ and }\mathcal{G}''=Z'\mathcal{F}''+Z''\mathcal{F}',\]
which can be written in matrix form as
\begin{equation}\label{sad_const}\left(\begin{array}{c}
\mathcal{G}'' \\
\mathcal{G}' \\
\end{array}\right)=\left(\begin{array}{cc}
Z'' & Z' \\
Z' & -Z'' \\
\end{array}\right)\left(\begin{array}{c}
\mathcal{F}' \\
\mathcal{F}'' \\
\end{array}\right).\end{equation}
Solving this relation for the imaginary parts of $\mathcal{F}$ and $\mathcal{G}$, we find that
\begin{equation}\label{min_const}\left(\begin{array}{c}
\mathcal{G}'' \\
\mathcal{F}'' \\
\end{array}\right)=\mathcal{L}\left(\begin{array}{c}
\mathcal{F}' \\
-\mathcal{G}' \\
\end{array}\right),\end{equation}
where
\[\mathcal{L}=\left(\begin{array}{cc}
Z''+Z'(Z'')^{-1}Z' & Z'(Z'')^{-1} \\
(Z'')^{-1}Z' & (Z'')^{-1} \\
\end{array}\right).\]
The matrix $\mathcal{L}$ is positive definite as long as $Z''$ is positive definite (see \cite{Milton_2009}). In terms of $L$ and $M$, this means that 
\[L''(x)>0\ \ \mbox{and}\ \ M''(x)>0\ \ \mbox{for all}\ x\in \Omega.\]
In practice it is only necessary that the values of $L$ and $M$ lie within any half plane of the complex plane.  This half plane can then be rotated so that it becomes the upper half plane (see Section~\ref{sad_princ}).

The approach in \cite{Richins_2012} was to use this constitutive relation and the corresponding energy functional 
\[\int_\Omega\left(\begin{array}{c}
\mathcal{F}' \\
-\mathcal{G}' \\
\end{array}\right)\cdot\left(\begin{array}{cc}
Z''+Z'(Z'')^{-1}Z' & Z'(Z'')^{-1} \\
(Z'')^{-1}Z' & (Z'')^{-1} \\
\end{array}\right)\left(\begin{array}{c}
\mathcal{F}' \\
-\mathcal{G}' \\
\end{array}\right)\ dx\]
to formulate a numerical method.  When this variational principle is discretized by the finite element method, the result is a system of equations that can be partitioned as a $3\times 3$ block matrix that consists of $N\times N$ blocks, where $N$ is the number of nodes in the computational grid.  A similar system of equations must be solved to find approximations for $u''$ and $v'$.  In all, to find $u'$ and $u''$, one must solve two positive definite systems of equations of size $3N\times 3N$. 

Here we develop a new method based on the saddle point variational principles in \cite{Milton_2009} that does not require that $v$ be solved for in order to find $u$, but is still based on solving positive definite systems of equations.  First, in Section~\ref{sad_princ}, we will analyze the saddle point variational principles from \cite{Milton_2009} upon which our method is based.  In Section~\ref{bdry_cond}, we will discuss the details of handling Dirichlet, Neumann, and Robin boundary conditions with these variational principles.  Section~\ref{err_bnd} contains the derivation of a standard bound on the error incurred when the Helmholtz equation is solved using a finite element method that discretizes the saddle point variational principle. Section~\ref{method} outlines the numerical method and discusses the conditioning of the system.  In Section~\ref{examples}, we provide several straightforward numerical explorations of the performance of the algorithm, as well as numerical verification of the error bound from Section~\ref{err_bnd}.  

\section{The Saddle Point Variational Principle}\label{sad_princ}

The derivation of the saddle point variational principle from \cite{Milton_2009} follows the same steps presented in the introduction for the minimization principle, the difference being that instead of continuing to the constitutive relation (\ref{min_const}), we stop at equation (\ref{sad_const}).  Assuming that $Z''$ is positive definite, from (\ref{sad_const}) we define the functional 
\begin{equation}Y(u',u'')=\int_\Omega\left(\begin{array}{c}
\mathcal{F}' \\
\mathcal{F}'' \\
\end{array}\right)\cdot\left(\begin{array}{cc}
Z'' & Z' \\
Z' & -Z'' \\
\end{array}\right)\left(\begin{array}{c}
\mathcal{F}' \\
\mathcal{F}'' \\
\end{array}\right)\ dx.\end{equation}
Let $u',u''\in H^1(\Omega)$ be the real and imaginary parts of a solution to the Helmholtz equation.  Let $s\in H_0^1(\Omega)$ and define 
\[\mathcal{S}=\left(\begin{array}{c}
\nabla s \\
s \\
\end{array}\right).\]
Then we have 
\[Y(u'+s,u'')=\int_\Omega\left(\begin{array}{c}
\mathcal{F}'+\mathcal{S} \\
\mathcal{F}'' \\
\end{array}\right)\cdot\left(\begin{array}{cc}
Z'' & Z' \\
Z' & -Z'' \\
\end{array}\right)\left(\begin{array}{c}
\mathcal{F}'+\mathcal{S} \\
\mathcal{F}'' \\
\end{array}\right)\ dx\]
\[=Y(u',u'')+2\int_\Omega\left(\begin{array}{c}
\mathcal{F}' \\
\mathcal{F}'' \\
\end{array}\right)\cdot\left(\begin{array}{cc}
Z'' & Z' \\
Z' & -Z'' \\
\end{array}\right)\left(\begin{array}{c}
\mathcal{S} \\
0 \\
\end{array}\right)\ dx+Y(s,0).\]
The integral in the line above can be rewritten as
\begin{equation}\int_\Omega\left(\begin{array}{c}
\mathcal{G}'' \\
\mathcal{G}' \\
\end{array}\right)\cdot\left(\begin{array}{c}
\mathcal{S} \\
0 \\
\end{array}\right)\ dx=\int_\Omega\left(\begin{array}{c}
-v' \\
-\nabla\cdot v' \\
\end{array}\right)\cdot\left(\begin{array}{c}
\nabla s \\
s \\
\end{array}\right)\ dx\end{equation}
\[=\int_\Omega\left[-v'\cdot\nabla s-\nabla\cdot v' s\right]\ dx=\int_\Omega-\nabla\cdot\left[v's\right]\ dx=\int_{\partial\Omega}-v'\cdot ns\ dS=0.\]
Therefore,
\[Y(u'+s,u'')=Y(u',u'')+\int_\Omega\mathcal{S}\cdot Z''\mathcal{S}\ dx,\]
and the last term must be nonnegative, since $Z''$ is assumed to be positive definite.  A similar calculation yields
\[Y(u',u''+s)=Y(u',u'')-\int_\Omega\mathcal{S}\cdot Z''\mathcal{S}\ dx.\]
This shows that $(u',u'')$ is at a saddle point of the functional $Y$.

Suppose that $(u',u'')$ is a saddle point of the functional $Y$.  Then the functional $Q(s',s'')=Y(u'+s',u''+s''),$ defined for all $s',s''\in H_0^1(\Omega),$ should have a saddle point at $s'=s''=0$.  A necessary condition for this to happen is that the first variation of $Q$ should vanish.  If
\[\mathcal{S}'=\left(\begin{array}{c}
\nabla s' \\
s' \\
\end{array}\right)\ \ \mbox{and}\ \ \mathcal{S}''=\left(\begin{array}{c}
\nabla s'' \\
s'' \\
\end{array}\right),\]
then we must have
\[0=\int_\Omega\left(\begin{array}{c}
\mathcal{F}' \\
\mathcal{F}'' \\
\end{array}\right)\cdot\left(\begin{array}{cc}
Z'' & Z' \\
Z' & -Z'' \\
\end{array}\right)\left(\begin{array}{c}
\mathcal{S}' \\
\mathcal{S}'' \\
\end{array}\right)\ dx.\]
After writing this equation out in terms of $u$, $s$, $L$, and $M$ and integrating by parts, we find that the integrals  
\[\int_\Omega(-\nabla\cdot L''\nabla u'+M''u'-\nabla\cdot L'\nabla u''+M'u'')s'\ dx\]
and
\[\int_\Omega(-\nabla\cdot L'\nabla u'+M'u'+\nabla\cdot L''\nabla u''-M'' u'')s''\ dx\]
must add to zero for any choices of $s'$ and $s''$ in $H_0^1(\Omega)$.  The real and imaginary parts of the equation $\nabla\cdot L\nabla u=Mu$ can be written as
\[\nabla\cdot L'\nabla u'-\nabla\cdot L''\nabla u''-M'u'+M''u''=0\]
and
\[ \nabla\cdot L'\nabla u''+\nabla\cdot L''\nabla u'-M'u''-M''u'=0.\]
Notice that the left-hand sides of these equations are just the opposites of the expressions multiplying $s'$ and $s''$ in the integrals above.  Since the result of the integral must be zero regardless of the choice of $s'$ and $s''$, the saddle point of $Y$ must be a solution to the Helmholtz equation.

So far, we have assumed that $Z''$ is positive definite, but it is often possible to use this method even when $L$ and $M$ do not have positive imaginary parts.  A solution of the equation 
\[\nabla\cdot L\nabla u=Mu\]
is also a solution to 
\begin{equation}\nabla\cdot e^{i\theta}L\nabla u=e^{i\theta}Mu,\end{equation}
where $\theta$ is a constant.  Therefore, to ensure that the imaginary part of $Z''$ is positive definite, we can apply a rotation so that the new coefficients $e^{i\theta}L$ and $e^{i\theta}M$ have positive imaginary parts.  The necessary conditions on $L$ and $M$ for the method to apply are that their values lie within one open half-plane.  That half-plane may then be rotated so that it becomes the upper half-plane.

\section{Boundary Conditions}\label{bdry_cond}
The calculations done above show that a saddle point of $Y$ satisfying $u'=f'$ and $u''=f''$ on $\partial\Omega$ is a solution of 
\[\left\{\begin{array}{ll}
\nabla\cdot L\nabla u=Mu & \mbox{in}\ \Omega \\
u=f & \mbox{on}\ \partial\Omega \\
\end{array}\right. .\]
We can also solve the Neumann problem
\begin{equation}\left\{\begin{array}{ll}
\nabla\cdot L\nabla u=Mu & \mbox{in}\ \Omega \\
v\cdot n=g & \mbox{on}\ \partial\Omega \\
\end{array}\right. .\end{equation}
Let $s',s''\in H^1(\Omega)$ be arbitrary test functions.  Then we have
\[0=\int_\Omega\left[(\nabla\cdot v'-\nabla\cdot v')s'+(-\nabla\cdot v''+\nabla\cdot v'')s''\right]\ dx\]
\[=\int_\Omega\left[-v'\cdot\nabla s'-\nabla\cdot v's'+v''\cdot\nabla s''+\nabla\cdot v''s''\right]\ dx+\int_{\partial\Omega}\left[s'v'\cdot n-s''v''\cdot n\right]\ dx\]
\[=\int_\Omega\left(\begin{array}{c}
\mathcal{G}'' \\
\mathcal{G}' \\
\end{array}\right)\cdot\left(\begin{array}{c}
\mathcal{S}' \\
\mathcal{S}'' \\
\end{array}\right)\ dx+\int_{\partial\Omega}\left[s'v'\cdot n-s''v''\cdot n\right]\ dS\]
\[=\int_\Omega\left(\begin{array}{c}
\mathcal{F}' \\
\mathcal{F}'' \\
\end{array}\right)\cdot\left(\begin{array}{cc}
Z'' & Z' \\
Z' & -Z'' \\
\end{array}\right)\left(\begin{array}{c}
\mathcal{S}' \\
\mathcal{S}'' \\
\end{array}\right)\ dx+\int_{\partial\Omega}\left[s'v'\cdot n-s''v''\cdot n\right]\ dS.\]
Therefore, in order to solve the Neumann problem, we solve the weak equation
\[\int_\Omega\left(\begin{array}{c}
\mathcal{F}' \\
\mathcal{F}'' \\
\end{array}\right)\cdot\left(\begin{array}{cc}
Z'' & Z' \\
Z' & -Z'' \\
\end{array}\right)\left(\begin{array}{c}
\mathcal{S}' \\
\mathcal{S}'' \\
\end{array}\right)\ dx=\int_{\partial\Omega}\left[-s'g'+s''g''\right]\ dS\ \ \mbox{for all}\ \ s',s''\in H^1(\Omega).\]

To solve the Robin problem
\begin{equation}\left\{\begin{array}{ll}
\nabla\cdot L\nabla u=Mu & \mbox{in}\ \Omega \\
u+a v\cdot n=g & \mbox{on}\ \partial\Omega \\
\end{array}\right.,\end{equation}
we begin with the weak form of the Neumann problem, which we will write as
\[0=\int_\Omega\left(\begin{array}{c}
\mathcal{F}' \\
\mathcal{F}'' \\
\end{array}\right)\cdot\left(\begin{array}{cc}
Z'' & Z' \\
Z' & -Z'' \\
\end{array}\right)\left(\begin{array}{c}
\mathcal{S}' \\
\mathcal{S}'' \\
\end{array}\right)\ dx+\int_{\partial\Omega}\left(\begin{array}{c}
s' \\
s'' \\
\end{array}\right)\cdot\left(\begin{array}{c}
v'\cdot n \\
-v''\cdot n \\
\end{array}\right)\ dS\]
We split the boundary condition into its real and imaginary parts as
\[u'+a'v'\cdot n-a''v''\cdot n=g'\]
\[u''+a'v''\cdot n+a''v'\cdot n=g'',\]
which we can write as
\[\left(\begin{array}{c}
u' \\
u'' \\
\end{array}\right)+\left(\begin{array}{cc}
a' & a'' \\
a'' & -a' \\
\end{array}\right)\left(\begin{array}{c}
v'\cdot n \\
-v''\cdot n \\
\end{array}\right)=\left(\begin{array}{c}
g' \\
g'' \\
\end{array}\right).\]
If the matrix in the equation above is called $W$, then
\[\left(\begin{array}{c}
v'\cdot n \\
-v''\cdot n \\
\end{array}\right)=-W^{-1}\left(\begin{array}{c}
u' \\
u'' \\
\end{array}\right)+W^{-1}\left(\begin{array}{c}
g' \\
g'' \\
\end{array}\right),\]
so the weak form of the equation with Robin boundary conditions is
\[\int_\Omega\left(\begin{array}{c}
\mathcal{F}' \\
\mathcal{F}'' \\
\end{array}\right)\cdot\left(\begin{array}{cc}
Z'' & Z' \\
Z' & -Z'' \\
\end{array}\right)\left(\begin{array}{c}
\mathcal{S}' \\
\mathcal{S}'' \\
\end{array}\right)\ dx-\int_{\partial\Omega}\left(\begin{array}{c}
u' \\
u'' \\
\end{array}\right)\cdot W^{-1}\left(\begin{array}{c}
s' \\
s'' \\
\end{array}\right)\ dS=\]
\[-\int_{\partial\Omega}\left(\begin{array}{c}
g' \\
g'' \\
\end{array}\right)\cdot W^{-1}\left(\begin{array}{c}
s' \\
s'' \\
\end{array}\right)\ dS.\]
The inverse of $W$ is
\begin{equation}W^{-1}=\frac{1}{-(a')^2-(a'')^2}\left(\begin{array}{cc}
-a' & -a'' \\
-a'' & a' \\
\end{array}\right)=\frac{1}{|a|^2}\left(\begin{array}{cc}
a' & a'' \\
a'' & -a' \\
\end{array}\right),\end{equation}
so if we require that $a'$ be negative, the matrix that results from discretizing the left-hand side will have the same block form as those that result from the other boundary conditions.  If $a'>0$, we can instead rotate so that $L''$ and $M''$ are negative so that the necessary block structure of the matrices is preserved.

Care must be taken with solving the Neumann and Robin problems when rotation is used, to ensure that the correct boundary conditions are enforced.  For example, if one desires to solve the Neumann problem
\[\left\{\begin{array}{ll}
\nabla\cdot L\nabla u=Mu & \mbox{in}\ \Omega \\
v\cdot n=g & \mbox{on}\ \partial\Omega \\
\end{array}\right.,\]
the rotated version of the problem is 
\[\left\{\begin{array}{ll}
\nabla\cdot e^{i\theta}L\nabla u=e^{i\theta}Mu & \mbox{in}\ \Omega \\
\tilde v\cdot n=e^{i\theta}g & \mbox{on}\ \partial\Omega \\
\end{array}\right.,\]
where $\tilde v=ie^{i\theta}L\nabla u$.

\section{Error Bound}\label{err_bnd}
We will make the following assumptions on $Z''$:
\begin{equation}\label{assume}
\begin{array}{l}\mbox{a. there is a constant }\gamma_1\mbox{ such that } |[Z'']_{ij}(x)|<\gamma_1\mbox{ for all }i,\mbox{ }j,\mbox{ and }x\in\Omega \\
\mbox{b. there is }\gamma_2\mbox{ such that }Z''(x)>\gamma_2 I\mbox{ for all }x\in\Omega.
\end{array}
\end{equation}
The requirements on $Z''$ are equivalent to requiring similar bounds on $L''$ and $M''$.  Define the space $V=[H^1(\Omega)]^2$, endowed with the norm
\begin{equation}\|(u',u'')\|_V=(\|u'\|_{H^1(\Omega)}^2+\|u''\|_{H^1(\Omega)}^2)^{\frac{1}{2}}.\end{equation}
Also, we will assume that $V_{N1}$ and $V_{N2}$ are finite dimensional subspaces of $H^1(\Omega)$, and that $V_N=V_{N1}\times V_{N2}$ is the space in which we seek our numerical solution.

Define a functional $f(s')$ for $s'\in H^1(\Omega)$ as
\[f(s')=\frac{1}{2}Y(s',u'')+Q(s',u'')=\frac{1}{2}\int_\Omega\left(\begin{array}{c}
\mathcal{S}' \\
\mathcal{F}'' \\
\end{array}\right)\cdot\left(\begin{array}{cc}
Z'' & Z' \\
Z' & -Z'' \\
\end{array}\right)\left(\begin{array}{c}
\mathcal{S}' \\
\mathcal{F}'' \\
\end{array}\right)\ dx+Q(s',u''),\]
where, in practice, $Q(s',u'')$ would contain terms that arise from the enforcement of boundary conditions and any inhomogeneous terms.  We will further divide the terms  as
\[f(s')=\frac{1}{2}B(s',s')-F(s',u''),\]
where 
\[B(s',s'')=\int_\Omega\mathcal{S}'\cdot Z''\mathcal{S}''\ dx\]
and $F(s',u'')$ contains the rest of the terms.  If $u'$ is a minimizer of $f(s')$, then $u'$ must satisfy the Euler-Lagrange equation
\begin{equation}B(u',s')=F(s',u'')\ \ \mbox{for all}\ \ s'\in H^1(\Omega).\end{equation}
Therefore, we can write
\[f(s')=\frac{1}{2}B(s',s')-F(s',u'')=B(u',u')-F(u',u'')+\frac{1}{2}B(s',s')-F(s',u'')\]
\[=B(u',u')-F(u',u'')+\frac{1}{2}B(s',s')-B(u',s')\]
\[=\frac{1}{2}B(u',u')-F(u',u'')+\frac{1}{2}B(u',u')-B(u',s')+\frac{1}{2}B(s',s')\]
\[=\frac{1}{2}B(u',u')-F(u',u'')+\frac{1}{2}B(u'-s',u'-s').\]

Suppose that $u_N'\in V_{N1}$ is such that 
\[f(u_N')=\min_{s\in V_{N1}}f(s').\]
Then
\begin{equation}\label{Bmin}B(u'-u'_N,u'-u'_N)^{\frac{1}{2}}=\min_{s'\in V_{N1}} B(u'-s',u'-s')^{\frac{1}{2}},\end{equation}
and the inequalities (\ref{assume}) imply that 
\[\sqrt{\gamma_2}\|s'\|_{H^1(\Omega)}\le\sqrt{B(s',s')}\le C\sqrt{\gamma_1}\|s'\|_{H^1(\Omega)}\ \ \mbox{for all}\ \ s'\in H^1(\Omega).\]
Applying these inequalities to both sides of (\ref{Bmin}) yields
\[\sqrt{\gamma_2}\|u'-u'_N\|_{H^1(\Omega)}\le \min_{s'\in V_{N1}} C\sqrt{\gamma_1}\|u'-s'\|_{H^1(\Omega)}.\]
Here and in what follows, $C$ will represent a constant that does not depend on $u'$, $u''$, or the grid spacing $h$.

In order to get the necessary bound, we must choose $s'$ properly.  Let $F_1$ be the orthogonal projection from $H^1(\Omega)$ onto $V_{N1}$.  Then $\|F_1\|_{B(H^1(\Omega),H^1(\Omega))}=1$, where $B(H^1(\Omega),H^1(\Omega))$ is the set of all bounded linear functions from $H^1(\Omega)$ to itself.  We then take $s'=F_1u'$ to obtain the inequality
\begin{equation}\label{u'bnd}\|u'-u'_N\|_{H^1(\Omega)}\le C\|u'-F_1u'\|_{H^1(\Omega)}.\end{equation}

If instead we use
\[f(s'')=\frac{1}{2}Y(u',s'')+Q(u',s'')=\frac{1}{2}\int_\Omega\left(\begin{array}{c}
\mathcal{F}' \\
\mathcal{S}'' \\
\end{array}\right)\cdot\left(\begin{array}{cc}
Z'' & Z' \\
Z' & -Z'' \\
\end{array}\right)\left(\begin{array}{c}
\mathcal{F}' \\
\mathcal{S}'' \\
\end{array}\right)\ dx+Q(u',s''),\]
and perform calculations similar to those above, we obtain the bound
\begin{equation}\label{u''bnd}\|u''-u''_N\|_{H^1(\Omega)}\le C\|u''-F_2u''\|_{H^1(\Omega)},\end{equation}
where $F_2$ is the orthogonal projection from $H^1(\Omega)$ onto $V_{N2}$.  

Combining inequalities (\ref{u'bnd}) and (\ref{u''bnd}), we find that
\[\|(u'-u_N',u''-u_N'')\|_V^2=\|u'-u_N'\|_{H^1(\Omega)}^2+\|u''-u_N''\|_{H^1(\Omega)}^2\]
\[\le C\left(\|u'-F_1u'\|_{H^1(\Omega)}^2+\|u''-F_2u''\|_{H^1(\Omega)}^2\right)=C\|(u'-F_1u',u''-F_2u'')\|_V^2\]
and consequently,
\[ \|(u',u'')-(u'_N,u''_N)\|_V\le C\|(u'-F_1u',u''-F_2u'')\|_V.\]

We partition $\Omega$ into subregions $e_l$, each of which can be viewed as a suitably shifted and rotated version of a reference element $\hat e$, so that there exist affine changes of variables $F_l(x)=B_lx+x_l$ such that $F_l(\hat e)=e_l$.  In what follows, a hat over a function will denote the corresponding function defined over the reference element $\hat e$ obtained by a change of variables.

We define the seminorm $|\cdot|_s$ by
\[|u|_s^2=[u,u]_s,\]
where
\begin{equation}\label{seminorm}[u,w]_s=\sum_{|\alpha|=s}\int_{\hat e} D^{\alpha}u\cdot D^{\alpha}w\ dx\end{equation}
and $\alpha$ is a multi-index.

 From \cite{Brezzi-Fortin} we get the inequality
\begin{equation}\label{inequalities}
c^{-1}h^{s-\frac{d}{2}}|w|_{s,e_l}\le|\hat w|_s\le ch^{s-\frac{d}{2}}|w|_{s,e_l}, \end{equation}
where $c$ is a constant, $w=\hat w\circ F_l^{-1}$, and the subscript $e_l$ denotes (\ref{seminorm}) with $e_l$ in place of $\hat e$.

We now recall the following lemma from \cite{Axelsson-Barker}:
\begin{lma}[Bramble-Hilbert Lemma] For some region $\Omega\subset\R^d$ and some integer $k\ge -1$, let there be given a bounded linear functional
\[f:H^{k+1}(\Omega)\rightarrow\R,\]
satisfying $|f(u)|\le\delta\|u\|_{H^{k+1}(\Omega)}$ for all $u\in H^{k+1}(\Omega)$ for some $\delta$ independent of $u$.  Suppose that $f(u)=0$ for all $u\in P_k(\bar\Omega)$.  Then there exists a constant $C$, dependent only on $\Omega$ such that 
\[|f(u)|\le C\delta|u|_{k+1},\ \ \ u\in H^{k+1}(\Omega).\]
\end{lma}

Let $s\in\{0,1\}$ and fix $w\in H^s(\hat e)$.  Define the functionals
\[L_1(\hat u)=[\hat u-F_1 \hat u,w]_s\ \mbox{and}\ L_2(\hat u)=[\hat u-F_2\hat u,w]_s.\]
Since
\[|L_j(\hat u)|\le|\hat u-F_j\hat u|_s|w|_s\le(|\hat u|_s+|F_j\hat u|_s)|w|_s\le(\|\hat u\|_{H^1(\hat e)}+\|F_j\hat u\|_{H^1(\hat e)})|w|_s\]
\[\le2\|\hat u\|_{H^1(\hat e)}|w|_s\le 2\|u\|_{H^{k+1}(\hat e)}|w|_s,\]
and $F_ju=u$ for polynomial functions $u$ in $V_{Nj}$ $(j=1,2)$, we see that the Bramble-Hilbert Lemma applies, and there exist constants such that 
\[|L_1(\hat u')|\le C|w|_s|\hat u'|_{k+1}\ \mbox{and}\ |L_2(\hat u'')|\le C|w|_s|\hat u''|_{k+1},\]
as long as $k$ is small enough so that all polynomials of degree less than or equal to $k$ are contained in the span of the basis functions representing $\hat u'$ and $\hat u''$.  Taking $w=\hat u'-F_1 \hat u'$ in the first inequality and $w=\hat u''-F_2\hat u''$ in the second yields
\begin{equation}|\hat u'-F_1\hat u'|_s\le C|\hat u'|_{k+1}\ \mbox{and}\ |\hat u''-F_2\hat u''|_s\le C|\hat u''|_{k+1}.\end{equation}

Assuming that $h\le 1$ and using inequality (\ref{inequalities}), we see that
\[|u'-F_1u'|_{s,e_l}\le Ch^{\frac{d}{2}-s}|\hat u'-F_1\hat u'|_s\le Ch^{\frac{d}{2}-s}|\hat u'|_{k+1}\le Ch^{k-s+1}|u'|_{k+1,e_l}\]
and
\[|u''-F_2u''|_{s,e_l}\le Ch^{\frac{d}{2}-s}|\hat u''-F_2 \hat u''|_s\le Ch^{\frac{d}{2}-s}|\hat u''|_{k+1}\le Ch^{k-s+1}|u''|_{k+1,e_l}\]

Consequently, the overall error satisfies
\[\|(u',u'')-(u'_N-u_N'')\|_V^2\le C\|(u'-F_1u',u''-F_2u'')\|_V^2\]
\[\le C\sum_l\left[|u'-F_1u'|_{0,e_l}^2+|u'-F_1u'|_{1,e_l}^2+|u''-F_2u''|_{0,e_l}^2+|u''-F_2u''|_{1,e_l}^2\right]\]
\[\le C\sum_l\left[h^{2k+2}|u'|_{k+1,e_l}^2+h^{2k}|u'|_{k+1,e_l}^2+h^{2k+2}|u''|_{k+1,e_l}^2+h^{2k}|u''|_{k+1,e_l}^2\right]\]
\[\le Ch^{2k}(|u'|_{k+1,\Omega}^2+|u''|_{k+1,\Omega}^2).\]

We have now proved
\begin{thm} Under the assumptions (\ref{assume}) on $Z''$, if the solution $(u',u'')\in [H^{k+1}(\Omega)]^2$ and the finite element subspace used in the numerical method contains $[P_k(\bar\Omega)]^2$, then there exists a constant $C$ such that the error satisfies
\[\|(u',u'')-(u'_N,u_N'')\|_V^2\le C h^{2k}(|u'|_{k+1,\Omega}^2+|u''|_{k+1,\Omega}^2),\]
where $h\le 1$ is the grid spacing.
\end{thm}

\section{The Numerical Method}\label{method}

To fix ideas, we will examine the numerical solution of the Dirichlet problem
\begin{equation}\left\{\begin{array}{ll}
\nabla\cdot L\nabla u=Mu & \mbox{in}\ \Omega \\
u=f & \mbox{on}\ \Omega \\
\end{array}\right.\end{equation}
The first step in solving the problem is to select a set of finite element basis functions. The numerical examples presented here will use a rectangular grid with bilinear basis functions. 

Regardless of how the basis is chosen, we will assume that the basis functions are labeled as $\{\psi_k\}$ and we assume that the solution has the form
\[\left(\begin{array}{c}
u' \\
u'' \\
\end{array}\right)=\left(\begin{array}{c}
\psi_0'+\sum\alpha_k'\psi_k \\
\psi_0''+\sum\alpha_k''\psi_k \\
\end{array}\right),\]
where $\psi_0'$ and $\psi_0''$ are auxiliary functions satisfying the boundary conditions $\psi_0'=f'$ and $\psi_0''=f''$ on $\partial\Omega$.  The weak form of the Euler-Lagrange equation for the saddle point variational principle is
\[0=\int_\Omega\left(\begin{array}{c}
\mathcal{F}' \\
\mathcal{F}'' \\
\end{array}\right)\cdot\left(\begin{array}{cc}
Z'' & Z' \\
Z' & -Z'' \\
\end{array}\right)\left(\begin{array}{c}
\mathcal{S}' \\
\mathcal{S}'' \\
\end{array}\right)\ dx\ \ \mbox{for all}\ \ s', s''\in H_0^1(\Omega),\]
where, as usual, 
\[\mathcal{S}'=\left(\begin{array}{c}
\nabla s' \\
s' \\
\end{array}\right)\ \ \mbox{and}\ \ \mathcal{S}''=\left(\begin{array}{c}
\nabla s'' \\
s'' \\
\end{array}\right).\]

We make the substitution above for $u'$ and $u''$ and let $s'$ and $s''$ be equal to each of the basis functions in turn.  In doing so, we arrive at a system of equations which has the block form
\begin{equation}\label{SaddleEqns}\left(\begin{array}{cc}
A_1 & A_2\\
A_2 & -A_1 \\
\end{array}\right)\left(\begin{array}{c}
\alpha' \\ 
\alpha'' \\
\end{array}\right)=\left(\begin{array}{c}
b_1 \\
b_2 \\
\end{array}\right),\end{equation}
where $A_1$ is positive definite.  The entries of the blocks in the coefficient matrix satisfy
\[[A_1]_{kj}=\int_\Omega\nabla\psi_k\cdot L''\nabla\psi_j\ dx+\int_\Omega\psi_k M''\psi_j\ dx\]
and
\[[A_2]_{kj}=\int_\Omega\nabla\psi_k\cdot L'\nabla\psi_j\ dx+\int_\Omega\psi_k M'\psi_j\ dx\]
The elements of the vector $b=(b_1,b_2)^T$ satisfy
\[[b_1]_j=-\int_\Omega\nabla\psi_0'\cdot L''\nabla\psi_j\ dx-\int_\Omega\psi_0'M''\psi_j\ dx-\int_\Omega\nabla\psi_0''\cdot L'\nabla\psi_j\ dx-\int_\Omega\psi_0''M'\psi_j\ dx\]
and
\[[b_2]_j=-\int_\Omega\nabla\psi_0'\cdot L'\nabla\psi_j\ dx-\int_\Omega\psi_0' M'\psi_j\ dx+\int_\Omega\nabla\psi_0''\cdot L''\nabla\psi_j\ dx+\int_\Omega\psi_0''M''\psi_j\ dx.\]

This system of equations (\ref{SaddleEqns}) is of saddle point type, and therefore there is a wide array of numerical methods that apply \cite{Golub_2005}.  Among the simplest is the following, based on Schur complements.  By using this approach, we reduce the problem from solving an indefinite $2N\times 2N$ system to solving two $N\times N$ positive definite systems.  We solve the second equation in (\ref{SaddleEqns}) for $\alpha''$ and substitute into the first to obtain
\begin{equation}\label{A1Primary}\begin{array}{l}
A_1\alpha''=-b_2+A_2\alpha' \\
(A_1+A_2A_1^{-1}A_2)\alpha'=b_1+A_2A_1^{-1}b_2 \\
\end{array}.\end{equation}
Because $A_1$ is positive definite and $A_2$ is symmetric, the coefficient matrices in both these systems of equations are positive definite. Equivalently, we can solve the second equation for $A_2$ and make the corresponding substitution into the first equation to obtain the system of equations
\begin{equation}\label{A2Primary}\begin{array}{l}
A_2\alpha'=b_2+A_1\alpha' \\
(A_2+A_1A_2^{-1}A_1)\alpha''=b_1-A_1A_2^{-1}b_2 \\
\end{array}.\end{equation}
The methods below can be adapted to this second system of equations under the assumption that $A_2$ is positive definite, which corresponds to $L$ and $M$ having positive real parts.  If the real parts of $L$ and $M$ are both positive, the problem can be rotated so that the imaginary parts become positive, so we will focus primarily on equations (\ref{A1Primary}).

While the matrix $A_1+A_2A_1^{-1}A_2$ is positive definite, it is also costly to store and to compute. For this reason, we use the preconditioned conjugate gradient (PCG) method to compute the solution to the system with this coefficient matrix, since this method only requires the ability to perform matrix-vector multiplication with the coefficient matrix.  As a preconditioner for  $A_1+A_2A_1^{-1}A_2$, we use the matrix $A_1$.  In this case, the preconditioned system has coefficient matrix
\[A_1^{-1}(A_1+A_2A_1^{-1}A_2)=I+(A_1^{-1}A_2)^2.\]
We can expect our system of equations to have the best conditioning when $A_1$ and $A_2$ are approximately the same, or alternatively when $\|A_1\|_2$ is much larger than $\|A_2\|_2$.

Systems with coefficient matrix $A_1$ appear explicitly in the algorithm, but must also be solved at each step when PCG is applied to the matrix $A_1+A_2A_1^{-1}A_2$, and there are many different ways in which this system can be solved.  In the numerical examples that follow, all the systems of equations of the form $A_1x=b$ are solved using PCG with an incomplete Cholesky factorization of $A_1$ as a preconditioner.  In essence, this introduces an inner and an outer PCG iteration in step~4 below. The following section illustrates how the total number of PCG iterations performed solving systems $A_1x=b$ is related to the size of the computational grid and the coefficients in the Helmholtz equation. 

The algorithm used here is as follows, though details such as the iterative solver or preconditioning method may be modified as desired:
\begin{enumerate}
\item Form the matrices $A_1$ and $A_2$.
\item Compute the right-hand side vectors $b_1$ and $b_2$.
\item Compute $w_1=b_1+A_2A_1^{-1}b_2$.
\item Solve $(A_1+A_2A_1^{-1}A_2)\alpha'=w_1$ using PCG with the preconditioner $A_1$.
\item Compute $w_2=-b_2+A_2\alpha'$.
\item Solve $A_1\alpha''=w_2$ by PCG with an incomplete Cholesky factorization of $A_1$ as preconditioner.
\end{enumerate}

This algorithm is completely implicit, and therefore is well suited for large-scale problems.  Because all that is required are sparse matrix-vector multiplications, parallel implementations of this algorithm can produce a significant speedup.  In particular, this algorithm could be implemented on a GPU cluster, where many graphics processing units (each of which contains many processing cores) are used in parallel to perform very fast computations.

In some situations, particularly those involving high frequency, $\|A_2\|_2$ is much larger than $\|A_1\|_2$, suggesting that we use formulation (\ref{A2Primary}).  However, $A_2$ is not positive definite and therefore neither is $A_2+A_1A_2^{-1}A_1$.  The basic algorithm outlined above can still be used in this case, provided that PCG is replaced by an iterative method that does not require positivity, such as GMRES.  As pointed out in \cite{Day_2001}, many solver packages are focused mainly on solving systems of equations with real matrices.  The approach above can be considered as an equivalent real formulation of the usual complex system of equations.

\subsection{Conditioning}\label{Cond} 
\begin{figure}[tbp]
   \centering
   \includegraphics[width=2.94in]{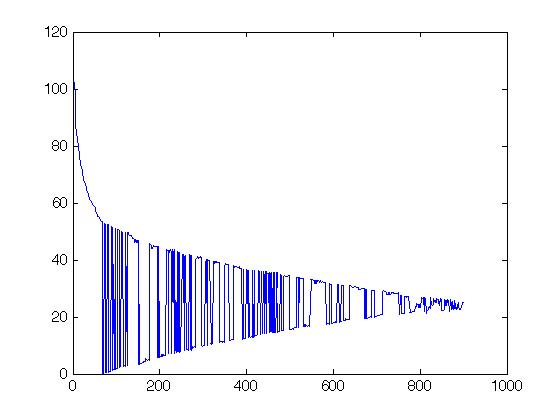} 
   \includegraphics[width=2.94in]{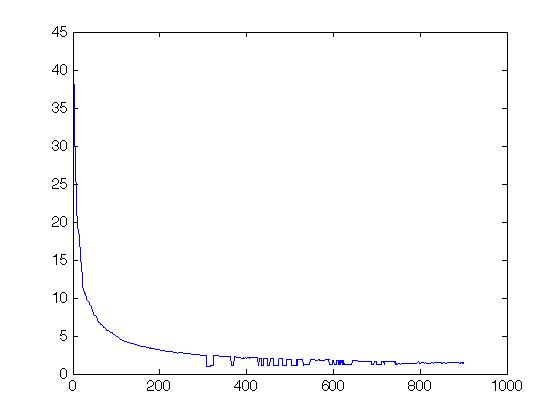}   
   \caption{The distribution of the eigenvalues of $A_1+A_2A_1^{-1}A_2$ (left) and the eigenvalues of $A_1^{-1}(A_1+A_2A_1^{-1}A_2)$ (right) for an example with $30\times 30$ computational grid.}
   \label{evals}
\end{figure}

In the numerical algorithm outlined above, we suggest that $A_1$ be used as a preconditioner for the system with matrix $A_1+A_2A_1^{-1}A_2$.  In Figure~\ref{evals}, we see the distribution of the eigenvalues of $A_1+A_2A_1^{-1}A_2$ and $A_1^{-1}(A_1+A_2A_1^{-1}A_2)$ for an example where the real and imaginary parts of $L$ and $M$ take on random values in the range $(0,10)$.  

Because the bulk of the work in this method comes from solving systems with matrix $A_1$, it is important that such systems can be effectively preconditioned.  A simple and effective choice is to use an incomplete Cholesky factorization of $A_1$ as the preconditioner, but there are many other preconditioning strategies that might be used.  If the algorithm is being implemented in parallel, a particularly useful strategy would be to use a sparse approximate inverse \cite{Benzi_1996}, which avoids the ``serial bottleneck" caused by having to perform back substitutions at each step in the PCG algorithm.

\begin{figure}[tbp]
\centering
\includegraphics[width=2.94in]{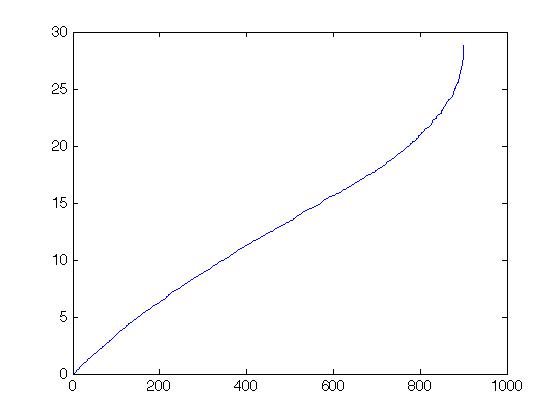}
\includegraphics[width=2.94in]{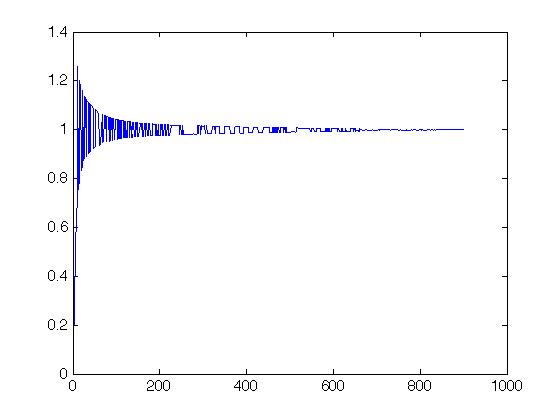}
\caption{The distribution of the eigenvalues of $A_1$ (left) and the eigenvalues of $(P^TP)^{-1}A_1$ (right), where $P$ is an incomplete Cholesky factorization of $A_1$ for an example with a $30\times 30$ computational grid.}
\label{A1evals}
\end{figure}

Figure~\ref{A1evals} shows the distribution of eigenvalues of $A_1$ before and after preconditioning.  The preconditioner used here is an incomplete Cholesky factorization of $A_1$ with drop tolerance $0.01$ and the real and imaginary parts of $L$ and $M$ take on random values is the range $(0,10)$.

\section{Numerical Results}\label{examples}

\begin{figure}[tbp]
\centering
\includegraphics[width=2.951in]{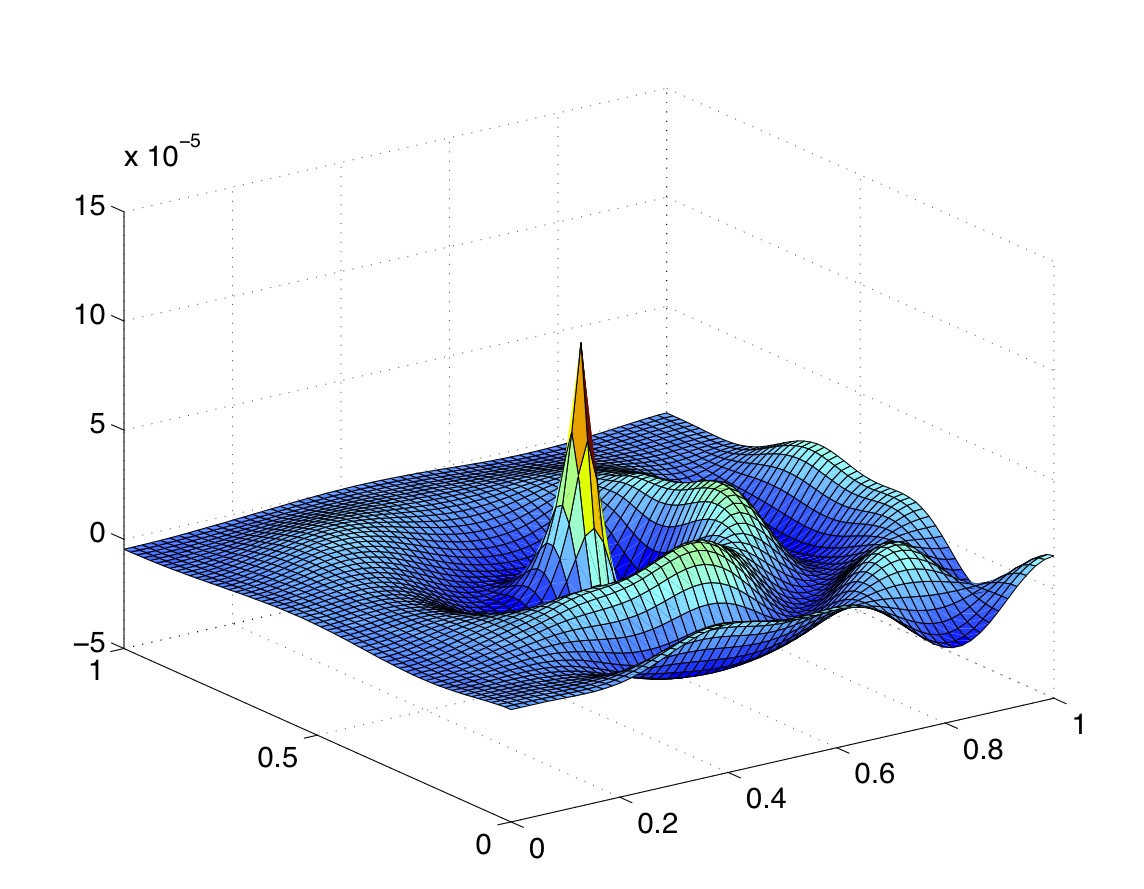}
\includegraphics[width=2.951in]{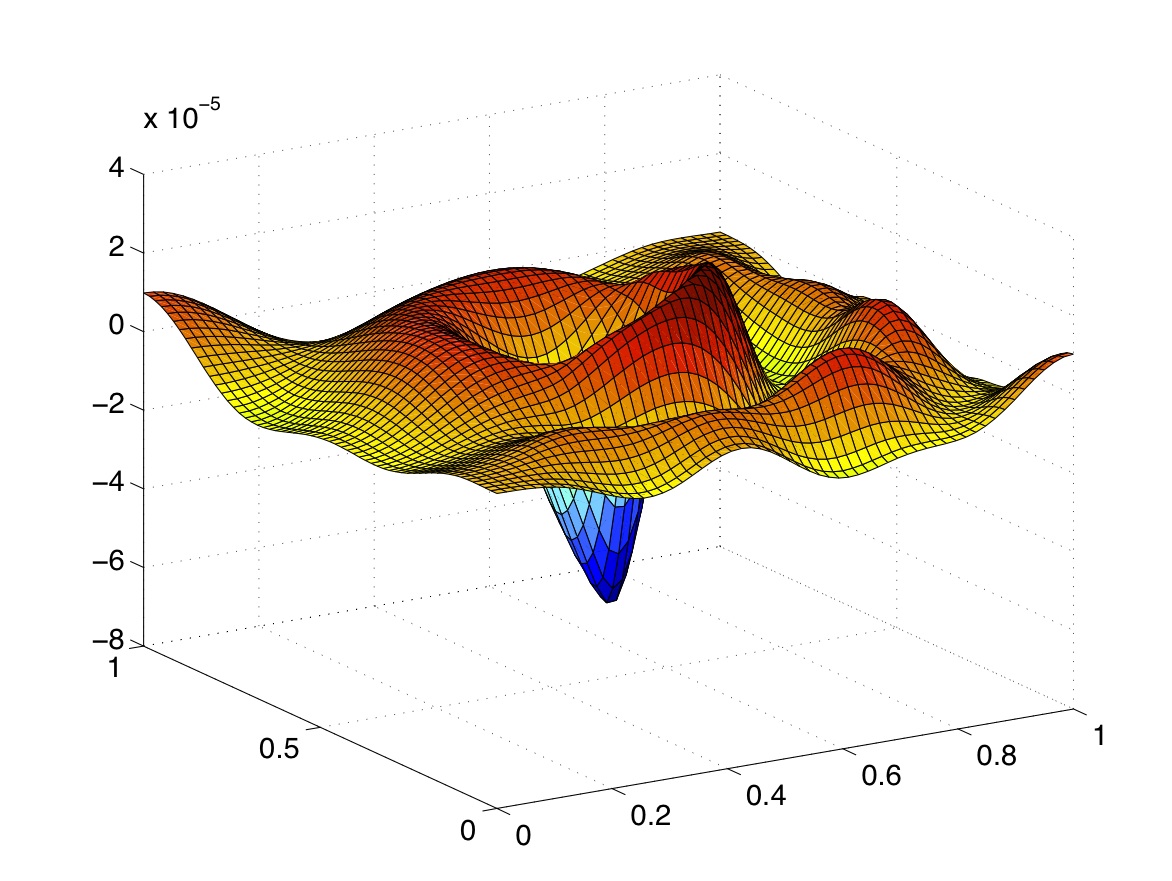}
\caption{The real (left) and imaginary (right) parts of a solution to the Helmholtz equation with a point source in the center of the domain and homogeneous Neumann boundary conditions.  The two phase material is chosen so that the frequency is higher below the line $y=x$ than above it.}
\label{DirEx}
\end{figure}

\begin{table}
\begin{center}
\begin{tabular}{|c|c|c|}
\hline
Grid & $h$ & $\|(u'-u'_N,u''-u''_N)\|_V^2$ \\
\hline
$32\times 32$ & 0.032258 & $1.4013\times 10^{-3}$ \\
$40\times 40$ & 0.025641 & $8.8176\times 10^{-4}$ \\
$50\times 50$ & 0.020408 & $5.5822\times 10^{-4}$ \\
$64\times 64$ & 0.015873 & $3.3785\times 10^{-4}$ \\
$70\times 70$ & 0.014493 & $2.8107\times 10^{-4}$ \\
$80\times 80$ & 0.012658 & $2.1502\times 10^{-4}$ \\
$90\times 90$ & 0.011236 & $1.6936\times 10^{-4}$ \\
$100\times 100$ & 0.010101 & $1.3663\times 10^{-4}$ \\
$128\times 128$ & 0.007874 & $8.3172\times 10^{-5}$ \\
$256\times 256$ & 0.003923 & $2.0663\times 10^{-5}$ \\
$512\times 512$ & 0.001957 & $5.1762\times 10^{-6}$ \\
\hline
\end{tabular}
\end{center}
\caption{The error in the finite element solution for various grid sizes.}
\label{error_table}
\end{table}

\begin{figure}[tbp]
\centering
\includegraphics[width=2.951in]{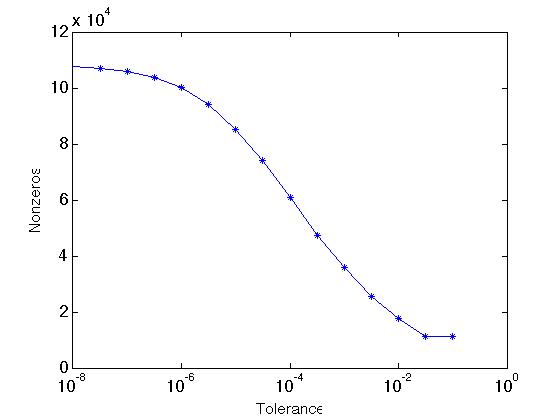}
\includegraphics[width=2.951in]{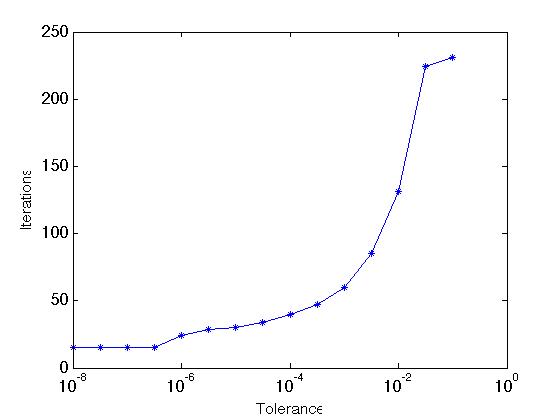}
\caption{Number of iterations and fill-in as a function of the drop tolerance used in the incomplete Cholesky factorization of $A_1$.}
\label{CholTest}
\end{figure}

\begin{figure}[tbp]
\centering
\includegraphics[width=4in]{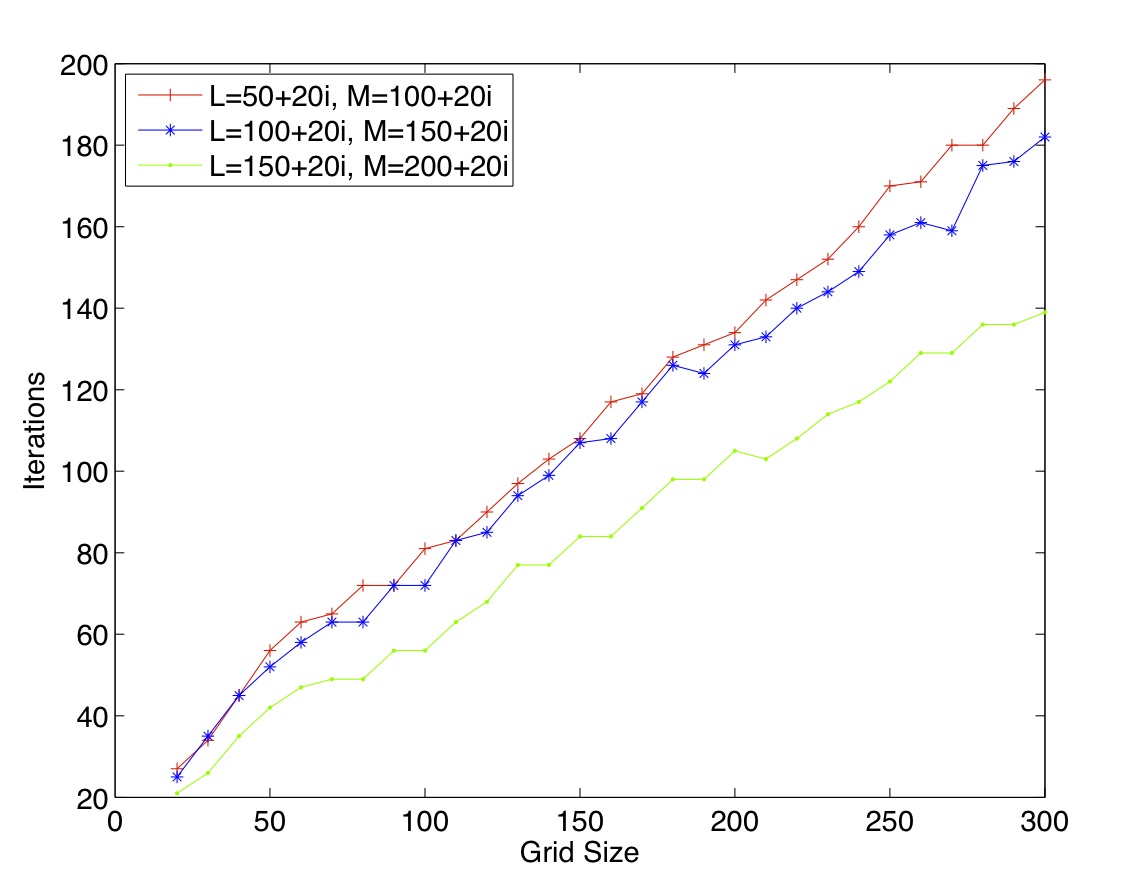}
\caption{The total number of PCG iterations required to solve the Helmholtz equation for several values of $L$  and $M$ as the size of the computational domain increases.}
\label{CompFig}
\end{figure}

\begin{figure}[tbp]
\centering
\includegraphics[width=4in]{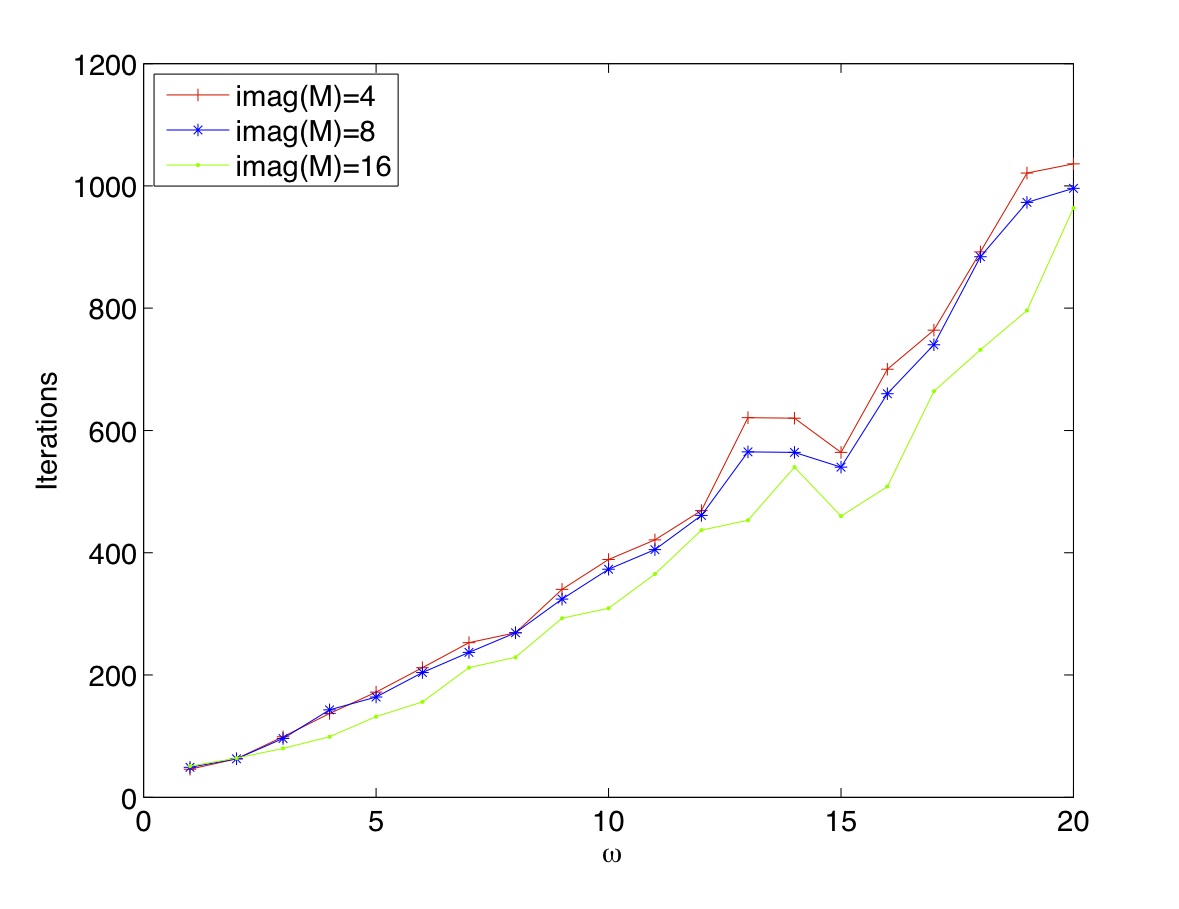}
\caption{The total number of PCG iterations required to solve the Helmholtz equation with $L=1$ and $M=-\omega^2+M''i$ for several values of $M''$.  The number of grid points per wavelength is held approximately constant at 10 as $\omega$ increases.}
\label{OmegaFig}
\end{figure}

\begin{figure}[tbp]
\centering
\includegraphics[width=2.75in]{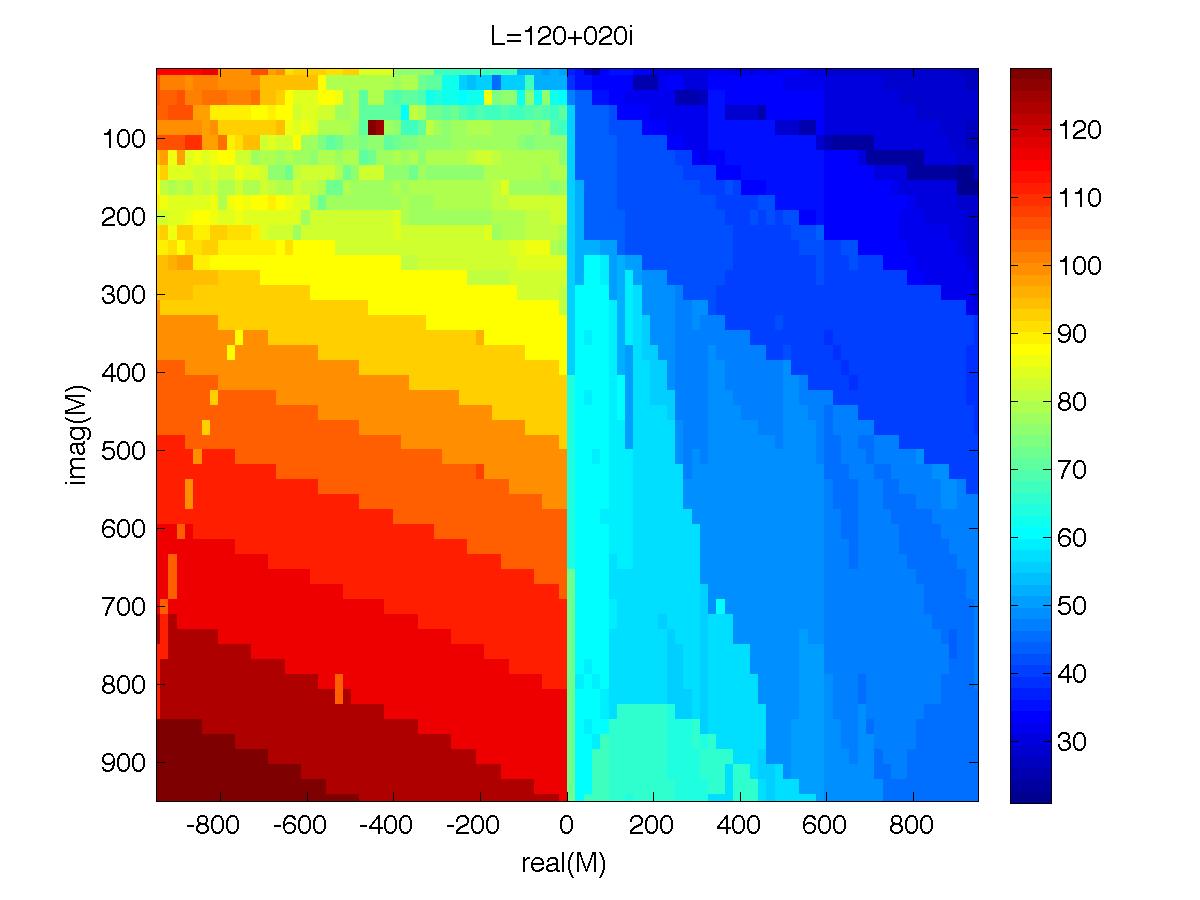}
\includegraphics[width=2.75in]{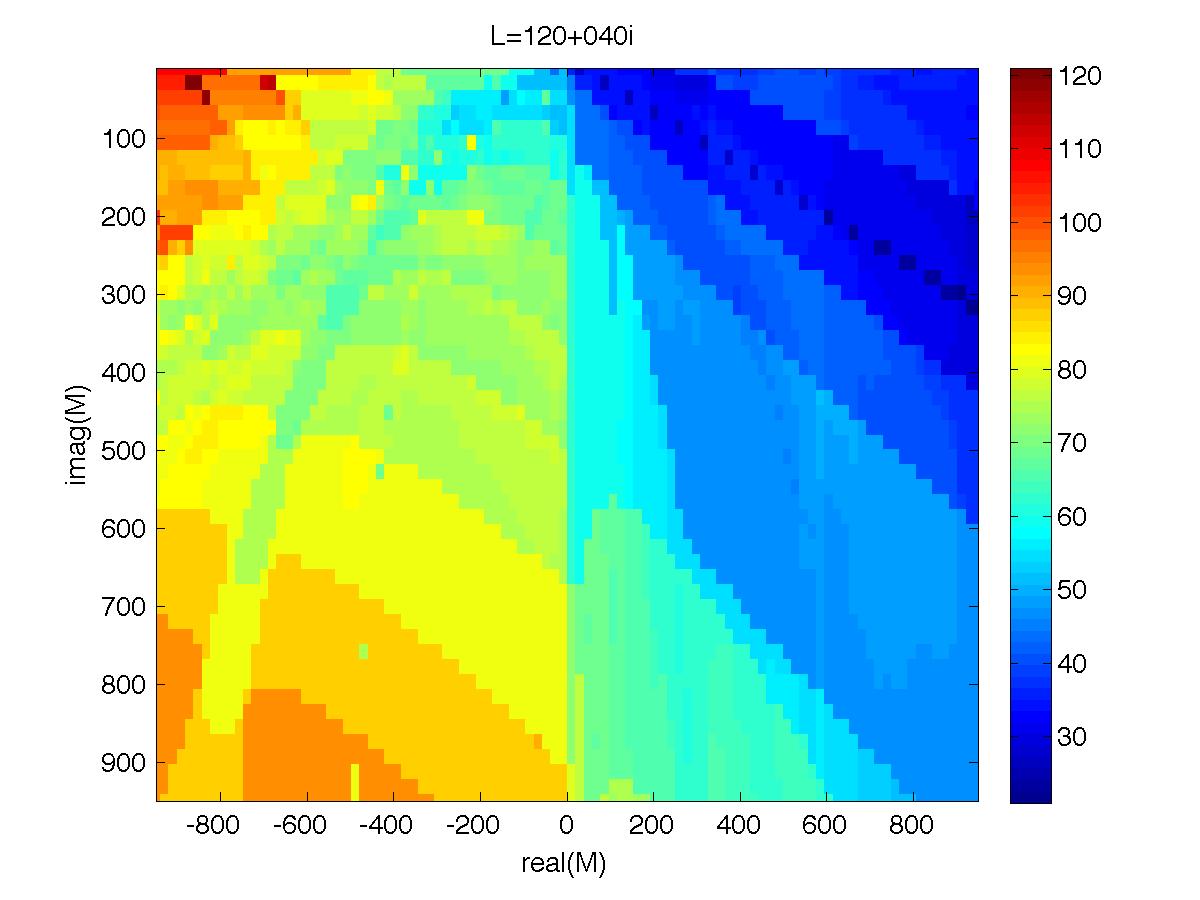}
\includegraphics[width=2.75in]{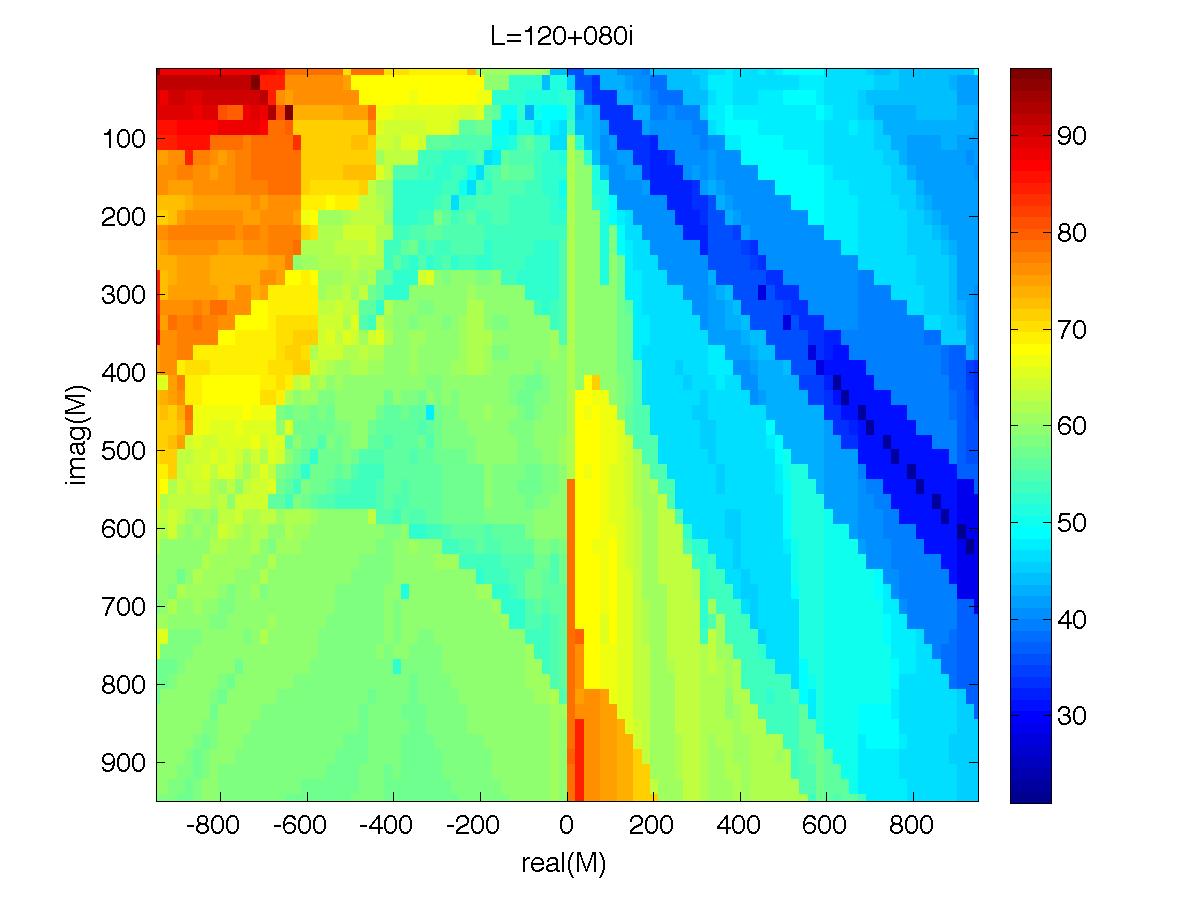}
\includegraphics[width=2.75in]{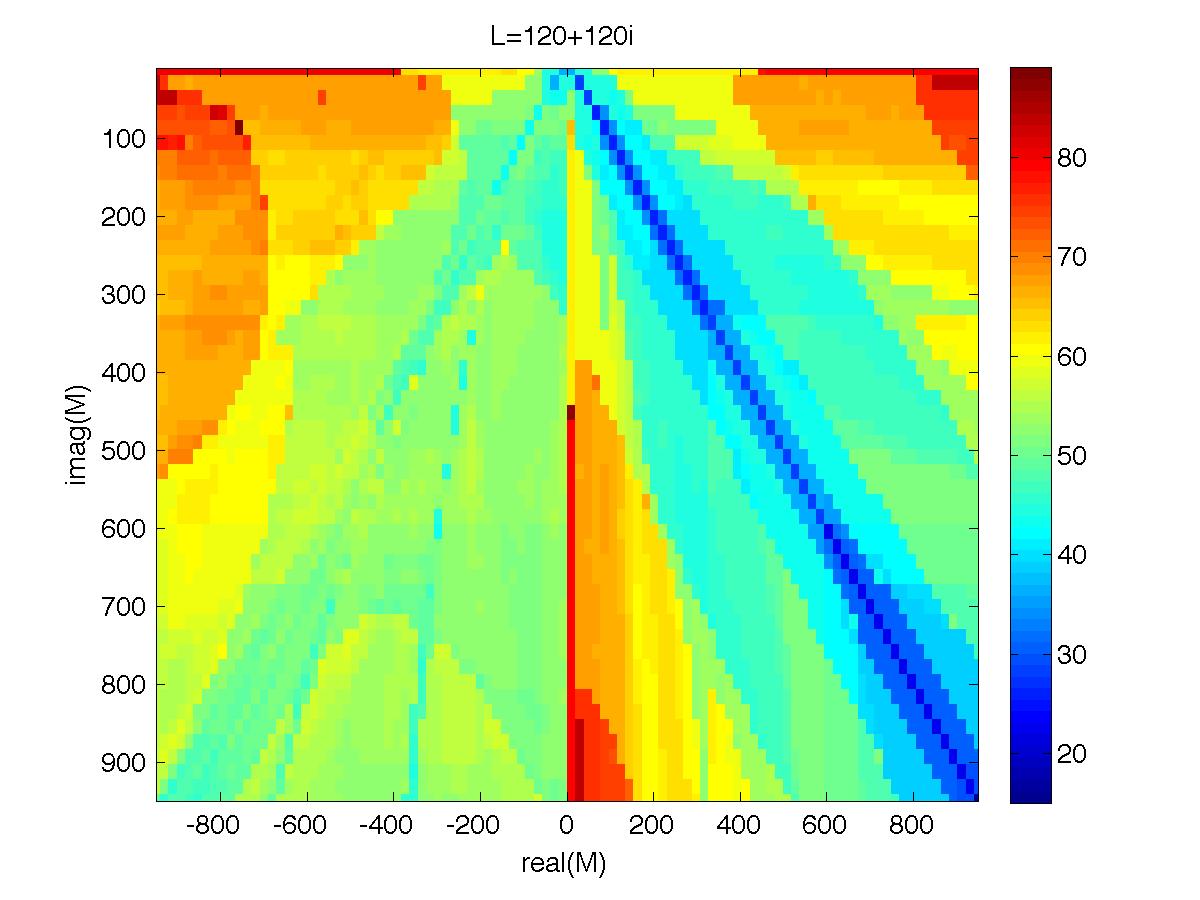}
\includegraphics[width=2.75in]{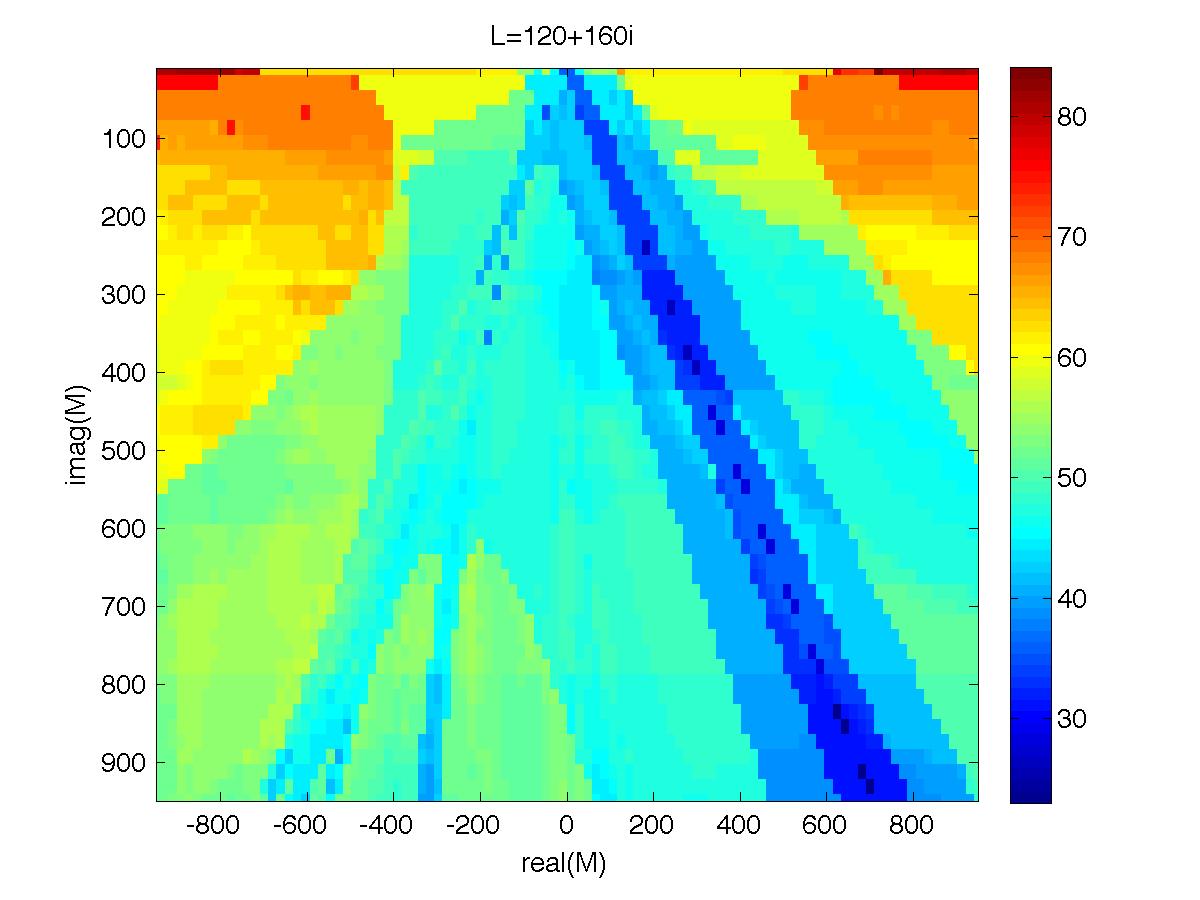}
\includegraphics[width=2.75in]{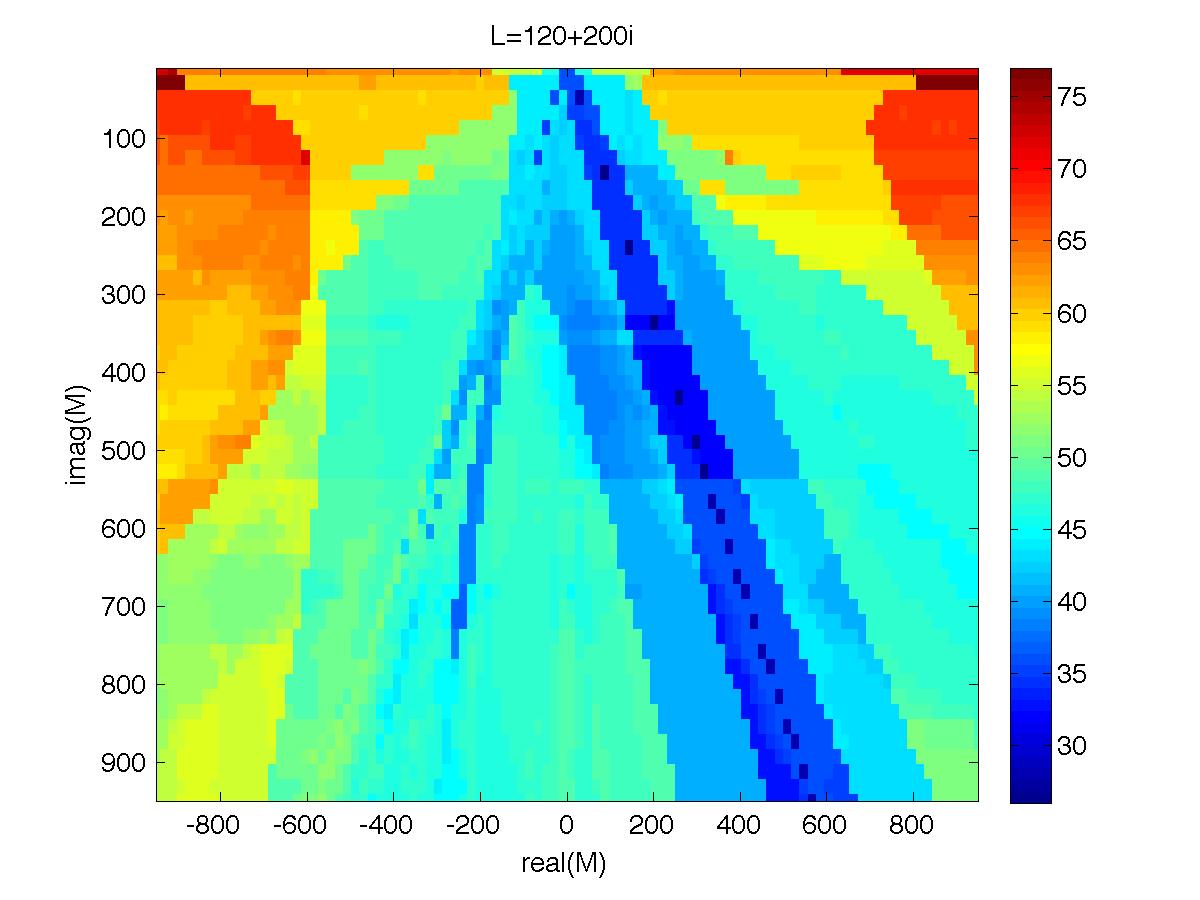}
\caption{The number of iterations required to solve the Helmholtz equation as the real and imaginary parts of $M$ vary for several values of $L$.  The grid size is fixed at $30\times 30$.}
\label{HeatFig}
\end{figure}

\begin{figure}[tbp]
\centering
\includegraphics[width=2.95in]{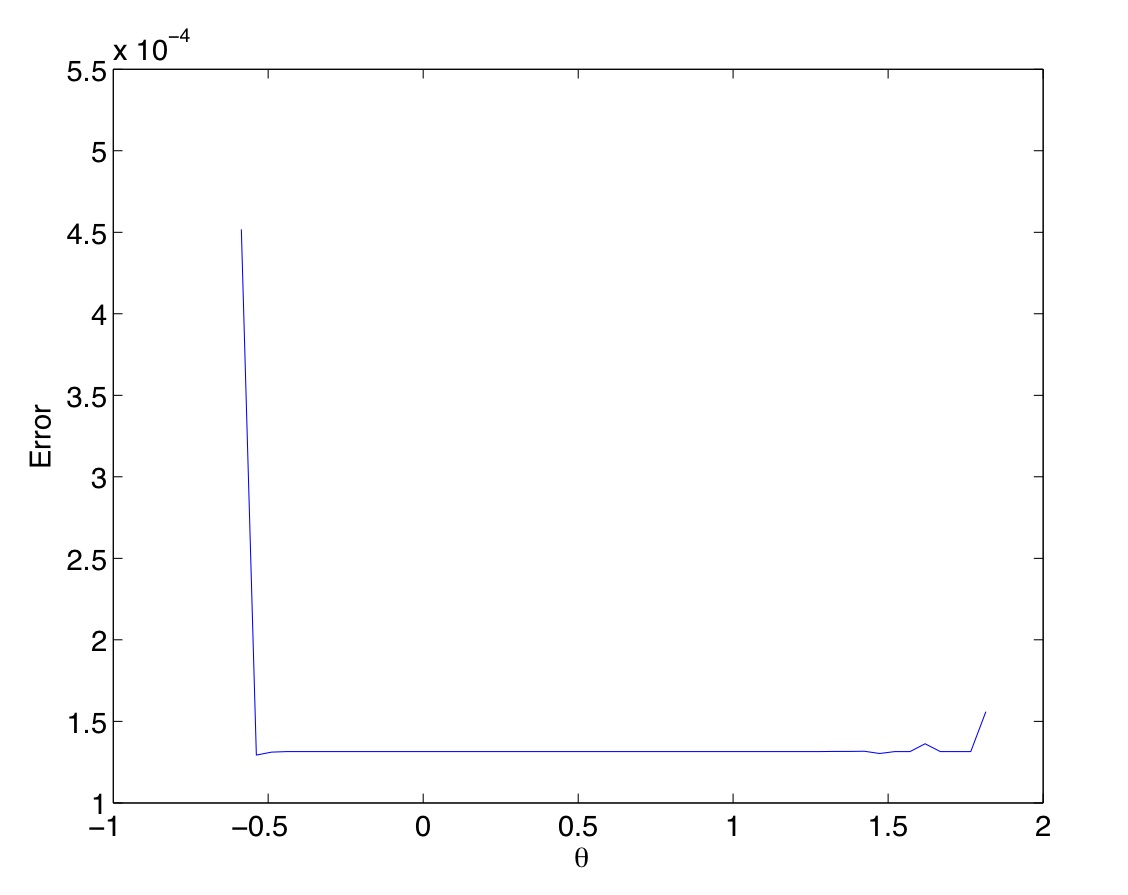}
\includegraphics[width=2.95in]{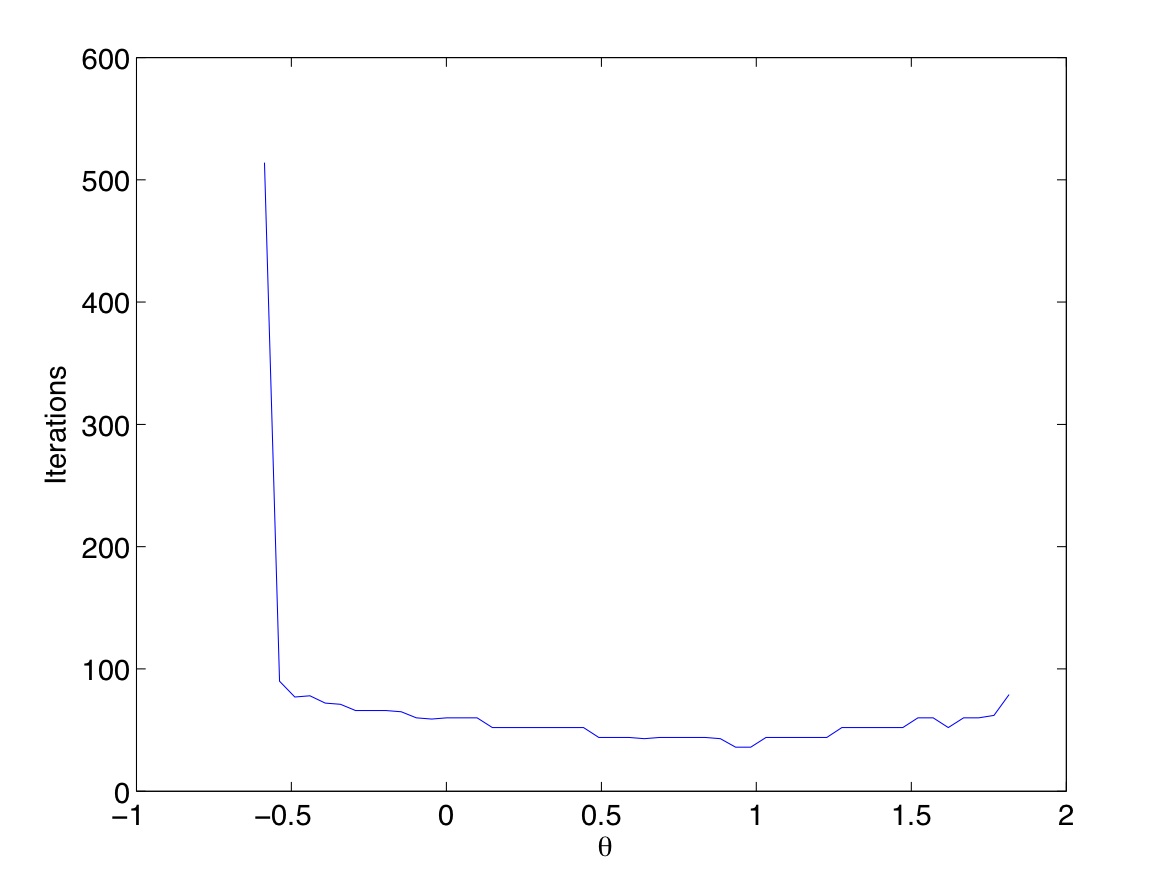}
\caption{The error and number of inner iterations in an example problem versus the value of $\theta$ used to rotate the problem.}
\label{rot_pic}
\end{figure}

In this section, we provide some demonstrations of the numerical solution of the Helmholtz problem
\[\nabla\cdot L\nabla u=Mu.\]

Figure~\ref{DirEx} demonstrates the application of this method to a non-homogeneous problem with variable coefficients.  The coefficients are chosen so that the frequency of the solution is higher in the lower right half of the unit square than in the upper left half, and the non-homogeneous term is a point source.  In the following numerical examples we will solve homogeneous problems with constant coefficients and attempt to quantify the convergence of the algorithm for different values of the coefficients.

Figure~\ref{CholTest} shows how the fill-in depends on the drop tolerance in the incomplete Cholesky factorization of $A_1$ (which is used in this section as preconditioner for systems with coefficient matrix $A_1$) for the problem with $L=1$ and $M=30-90i$ (before rotation) and the number of PCG iterations necessary to solve the same problem with a tolerance of $1\times 10^{-4}$ on the relative residual.

In Figure~\ref{CompFig}, we see the total number of PCG iterations necessary to solve all the systems with coefficient matrix $A_1$ for several different values of the coefficients in the problem as the size of the computational grid increases.  The tolerances for the PCG algorithm is $1\times 10^{-6}$ and the drop tolerance for the incomplete Cholesky factorization of $A_1$ is $1\times 10^{-4}$.  It should be noted that in this case the growth in iterations happens entirely within the inner PCG iterations.  The number of outer iterations required was either 2 or 3 in every instance.

Figure~\ref{OmegaFig} shows how the number of iterations is related to the frequency for problems where $M$ is in the left-hand side of the upper half plane.  In this situation, we cannot take advantage of the fact that $\|A_2\|_2>\|A_1\|_2$ if we wish to solve only positive definite systems because $A_2$ is not positive definite.  The number of grid points per wavelength is held approximately constant at 10 points per wavelength as $\omega$ grows.

The graphs in Figure~\ref{HeatFig} show in more detail how the number of iterations required to solve the Helmholtz equation depend on the coefficients in the problem.  In order to get maximum advantage from the preconditioning strategy outlined in Section~\ref{Cond}, when both coefficients $L$ and $M$ are in the first quadrant we choose formulation (\ref{A1Primary}) when $\|A_1\|_\infty>\|A_2\|_\infty$, and we choose formulation (\ref{A2Primary}) when the reverse inequality holds.  Unfortunately, we must use formulation (\ref{A1Primary}) when $L$ is in the first quadrant and $M$ is in the second quadrant because in this situation only $A_1$ is positive definite.

Table~\ref{error_table} shows the relationship between the error and the grid spacing in a problem with $L=-0.25+0.25i$ and $M=0.1+0.3i$ and Dirichlet boundary conditions.  In Figure~\ref{rot_pic}, the result of rotation on an example with $L=3+2i$ and $M=1+4i$ is shown.  The error and number of iterations remain nearly constant until $\theta$ is such that one of the imaginary parts of the rotated coefficients approaches zero.

\section{Conclusion}

By formulating a finite element method through the saddle-point variational principles of Milton, Seppecher, and Bouchitt\'{e}, we are able to solve boundary value problems for the complex Helmholtz equation by solving symmetric positive definite systems of equations.  The method is based on using elimination on the block structure of the finite element matrix to produce two smaller systems of equations, both of which have positive definite coefficient matrices.  The systems can then be solved using purely iterative methods. This method applies to a large class of problems, especially in light of the ability to ``rotate" the coefficients of a given problem to fit the assumptions of the algorithm.    

It should be emphasized that the method developed here does not only apply to the Helmholtz equation.  In \cite{Milton_2009}, there are similar variational principles given for the time-harmonic Maxwell equations and the equations of linear elasticity in lossy materials.  The ideas presented here can easily be adapted to these situations.  Also, the original variational principles of this type, developed by Cherkaev and Gibianski in \cite{Cherkaev_1994}, can be used to apply this numerical method to the complex Poisson equation.

As with the previous minimization-based method, the variational principles upon which this method is based remain valid as long as $L$ and $M$ have positive imaginary part, but the conditioning of the system deteriorates and the error incurred increases as $L$ and $M$ come close to violating this condition.

There is still more study necessary to determine the conditions under which this approach is competitive with other methods already in use.  Also, it is worthwhile to consider other boundary conditions in addition to the ones presented herein, such as a PML \cite{Turkel_2000}.  Also, the application of this method to problems with a non-local boundary condition, such as those considered in \cite{Bao_2005} may also be explored.

In Section~\ref{examples}, the preconditioning method used in the inner iterations was simply an incomplete Cholesky factorization.  To the extent that the growth in iterations in the inner iterations can be controlled, this method will become more attractive.  Future work in this direction will be to compare potential preconditioning methods and their performance in the overall algorithm, including multigrid, sparse approximate inverse \cite{Benzi_1996}, and sweeping preconditioners \cite{Engquist_2011}, and also to compare the amount of work required when this method is implemented in parallel to more standard methods of solving Helmholtz equations.

\section{Acknowledgements}

The author would like to thank David Dobson for his helpful comments during the preparation of the manuscript and Gang Bao for emphasizing the need for the improvements contained in this paper.

\bibliographystyle{amsplain}
\bibliography{paperbibliography}

\end{document}